\documentclass[12pt]{article}
\usepackage{makeidx}
\usepackage{amssymb}
\usepackage{amsmath}
\usepackage{amsfonts}
\usepackage{mitpress}

\newdimen\dummy
\dummy=\oddsidemargin
\input{tcilatex}

\begin{document}

\title{Analyse alg\'{e}brique appliqu\'{e}e \`{a} l'automatique}
\author{Henri Bourl\`{e}s}
\maketitle

\begin{abstract}
The expression "Algebraic Analysis" was coined by Mikio Sato. It consists of
using algebraic notions to solve analytic problem. The origin of Algebraic
Analysis is Algebraic Geometry as was developed by Alexander Grothendieck
and his school. Mimicking the introduction of Grothendieck's EGA (changing
only a few words) one obtains a good definition of the modern theory of
linear systems, as developed by Michel Fliess, Ian Willems, Ulrich Oberst
and others.
\end{abstract}

\sloppy

\section{Introduction}

L'expression "Analyse Alg\'{e}brique" est due \`{a} M. Sato, l'inventeur des
hyperfonctions. On pourrait dire que qu'appartiennent \`{a} l'Analyse Alg%
\'{e}brique tous les outils math\'{e}matiques utilisant de l'Alg\`{e}bre
pour \'{e}tudier des objets relevant de l'Analyse. Cette d\'{e}finition est
bien entendu bien trop g\'{e}n\'{e}rale. Si le "p\`{e}re" de l'analyse alg%
\'{e}brique est M. Sato, son "grand-p\`{e}re" est A. Grothendieck,
l'inventeur de la G\'{e}om\'{e}trie Alg\'{e}brique moderne. C'est ce que
nous allons montrer bri\`{e}vement dans ce qui suit.

\section{Le paradigme de la G\'{e}om\'{e}trie Alg\'{e}brique}

\subsection{La G\'{e}om\'{e}trie Alg\'{e}brique classique}

Classiquement, la G\'{e}om\'{e}trie Alg\'{e}brique s'int\'{e}resse aux
"ensembles alg\'{e}briques" d\'{e}finis comme suit: soit $\mathbf{k}$ un
corps commutatif et $\mathbf{k}\left[ T\right] $ l'anneau des polyn\^{o}mes
par rapport \`{a} une famille d'ind\'{e}termin\'{e}es $T=\left(
T_{1},...,T_{k}\right) ,$ \`{a} coefficients dans $\mathbf{k}.$ \ Soit
d'autre part $\mathbf{A}$ une $\mathbf{k}-$alg\`{e}bre et $\left\{
1,...,q\right\} $ un ensemble d'indices. \ Un sous-ensemble alg\'{e}brique $%
S $ de $\mathbf{A}^{k}$ est un ensemble%
\begin{equation}
V_{S}\left( \mathbf{A}\right) =\left\{ a=\left( a_{1},...,a_{k}\right) \in 
\mathbf{A}^{k}:F_{j}\left( a\right) =0,\forall j\in \left\{ 1,...,q\right\}
\right\}  \label{eq-ens-alg}
\end{equation}%
o\`{u} $F_{j}\left( T\right) \in \mathbf{P}_{k}\triangleq \mathbf{k}\left[ T%
\right] $ pour tout $j\in \left\{ 1,...,q\right\} $. \ Autrement dit, un
sous-ensemble de $\mathbf{A}^{k}$ est alg\'{e}brique quand il est d\'{e}fini
par des \'{e}quations polynomiales.

Ceci peut \^{e}tre g\'{e}n\'{e}ralis\'{e} en rempla\c{c}ant les ensembles
finis d'indices $\left\{ 1,...,k\right\} $ et $\left\{ 1,...,q\right\} $ par
des ensembles infinis.

\subsection{La r\'{e}volution grothendieckienne\label{par-revol-groth}}

On montre facilement que%
\begin{equation*}
V_{S}:\mathbf{A}\mapsto V_{S}\left( \mathbf{A}\right)
\end{equation*}%
est un foncteur covariant de la cat\'{e}gorie des $\mathbf{k}$-alg\`{e}bres
dans celle des ensembles. On peut distinguer deux \'{e}tapes dans l'\'{e}%
tude de l'ensemble alg\'{e}brique $V_{S}\left( \mathbf{A}\right) :$

\begin{itemize}
\item[(i)] l'\'{e}tude du foncteur $V_{S}.$

\item[(ii)] l'\'{e}tude de "l'immersion affine" particuli\`{e}re%
\begin{equation*}
V_{S}\left( \mathbf{A}\right) \rightarrow \mathbf{A}^{k}.
\end{equation*}
\end{itemize}

Remarquons maintenant que $V_{S}\left( \mathbf{A}\right) $ est l'ensemble
des points de $\mathbf{A}^{k}$ qui s'annulent sur l'id\'{e}al $\mathfrak{I}$
de $\mathbf{P}_{k}$ engendr\'{e} par les $F_{j}\left( T\right) .$ \ Il est
maintenant ais\'{e} de montrer qu'il existe un isomorphisme fonctoriel%
\begin{equation*}
V_{S}\left( \mathbf{A}\right) \cong \limfunc{Hom}\nolimits_{\mathbf{k}\text{%
-alg}}\left( \mathbf{P}_{k}/\mathfrak{I},\mathbf{A}\right) .
\end{equation*}

Comme toute $\mathbf{k}$-alg\`{e}bre commutative de pr\'{e}sentation finie $%
\mathbf{K}$ est de la forme $\mathbf{P}_{k}/\mathfrak{I},$ o\`{u} $\mathfrak{%
I}$ est un id\'{e}al de type fini de $\mathbf{P}_{k}$ (\cite{Bourbaki-Al-I}, 
\S III.2, n$%
%TCIMACRO{\U{b0}}%
%BeginExpansion
{{}^\circ}%
%EndExpansion
$8 et 9),\ les $V_{S}$ s'identifient aux foncteurs repr\'{e}sentables%
\begin{equation*}
V_{\mathbf{K}}:\mathbf{A}\mapsto \limfunc{Hom}\nolimits_{\mathbf{k}\text{-alg%
}}\left( \mathbf{K},\mathbf{A}\right) .
\end{equation*}

Le foncteur $\mathbf{K}\mapsto V_{\mathbf{K}}$ est contravariant et la cat%
\'{e}gorie des foncteurs%
\begin{equation*}
\mathbf{k}\text{-\textbf{alg}}^{cfp}\rightarrow \mathbf{Ens}
\end{equation*}%
associ\'{e}s aux ensembles $\left( \ref{eq-ens-alg}\right) $ est \'{e}%
quivalente \`{a} la cat\'{e}gorie oppos\'{e}e de la cat\'{e}gorie des $%
\mathbf{k}$-alg\`{e}bres commutatives de pr\'{e}sentation finie $\mathbf{k}$-%
\textbf{alg}$^{cfp}$ en associant \`{a} toute $\mathbf{k}$-alg\`{e}bre $%
\mathbf{K}$ (commutative et de pr\'{e}sentation finie) le foncteur $V_{%
\mathbf{K}}.$

Comme on l'a mentionn\'{e} plus haut, on peut consid\'{e}rer des ensembles
d'indices infinis, et on parvient donc au point de vue suivant \cite%
{Grothendieck-Dieudonne}:

\emph{Le but initial de la G\'{e}om\'{e}trie Alg\'{e}brique moderne \'{e}%
quivaut \`{a} l'\'{e}tude des }$\mathbf{k}$\emph{-alg\`{e}bres commutatives }%
$\mathbf{K}$\emph{.}

Telle est "l'alg\'{e}brisation" de la G\'{e}om\'{e}trie Alg\'{e}brique qui a
servi de point de d\'{e}part aux travaux d'A. Grothendieck et de son \'{e}%
cole (th\'{e}orie des sch\'{e}mas, etc.).

\section{Th\'{e}orie des syst\`{e}mes lin\'{e}aires}

\subsection{Point de vue classique\label{par-syst-classiques}}

Consid\'{e}rons pour fixer les id\'{e}es un syst\`{e}me lin\'{e}aire $S$ de
dimension finie. Un tel syst\`{e}me peut \^{e}tre donn\'{e} par une repr\'{e}%
sentation d'\'{e}tat (notion due \`{a} Kalman \cite{Kalman-1960})%
\begin{eqnarray*}
\dot{x} &=&A\left( t\right) x+B\left( t\right) u \\
y &=&C\left( t\right) x+D\left( t\right) u
\end{eqnarray*}%
o\`{u} $A,B,C$ et $D$ sont des matrices qu'on supposera, par exemple, \`{a}
coefficients analytiques, et o\`{u} $x,u$ et $y$ sont, respectivement, l'%
\'{e}tat, la commande et la sortie du syst\`{e}me.

La donn\'{e}e d'un syst\`{e}me par une repr\'{e}sentation d'\'{e}tat n'est
toutefois pas toujours naturelle. Un type de repr\'{e}sentation plus g\'{e}n%
\'{e}ral est d\^{u} \`{a} Rosenbrock \cite{Rosenbrock-book}:%
\begin{eqnarray*}
D\left( \partial \right) \xi &=&N\left( \partial \right) u, \\
y &=&Q\left( \partial \right) \xi +W\left( \partial \right) u
\end{eqnarray*}%
o\`{u} $D,N,Q$ et $W$ sont des matrices polynomiales par rapport \`{a} $%
\partial \triangleq d/dt$, les polyn\^{o}mes \'{e}tant suppos\'{e}s \`{a}
coefficients analytiques pour rester dans le m\^{e}me contexte que
ci-dessus, et o\`{u} $\xi ,u$ et $y$ sont, respectivement, l'\'{e}tat
partiel, la commande et la sortie du syst\`{e}me.

Dans certains cas, on ne sait pas \emph{a priori} quelles sont les variables
du syst\`{e}me qui sont appropri\'{e}es pour constituer la commande ou la
sortie (ce point a \'{e}t\'{e} soulign\'{e} par Willems \cite{Willems 91}).
Il convient donc de consid\'{e}rer une repr\'{e}sentation encore plus g\'{e}n%
\'{e}rale, de la forme%
\begin{equation}
R\left( \partial \right) w=0  \label{syst-Willems}
\end{equation}%
o\`{u} $R$ est une matrice polynomiale du m\^{e}me type que ci-dessus, de
dimension $q\times k$ par exemple.

Il convient de fixer dans quel espace $W$ la variable $w$ va pouvoir varier
(espace des fonctions analytiques, ou des fonctions ind\'{e}finiment d\'{e}%
rivables, ou des distributions, etc.) et Willems d\'{e}finit le "behavior"
(terme que l'on pourrait traduite en fran\c{c}ais par "comportement", mais
non sans ambigu\"{\i}t\'{e}) du syst\`{e}me $S$ comme suit:%
\begin{equation*}
\mathfrak{B}_{S}\left( W\right) =\left\{ \mathbf{w}\in W^{k}:R\left(
\partial \right) \mathbf{w}=0\right\} .
\end{equation*}

\subsection{Le point de vue de l'Analyse Alg\'{e}brique}

Pour plus de g\'{e}n\'{e}ralit\'{e}, r\'{e}\'{e}crivons $\left( \ref%
{syst-Willems}\right) $ sous la forme%
\begin{equation}
Rw=0  \label{def-matricielle-systeme}
\end{equation}%
o\`{u} $R$ est une matrice de dimension $q\times k$ \`{a} coefficients dans
un anneau d'op\'{e}rateurs $\mathbf{D.}$ On supposera que $\mathbf{D}$ est
une $%
%TCIMACRO{\U{211d} }%
%BeginExpansion
\mathbb{R}
%EndExpansion
$-alg\`{e}bre, \'{e}ventuellement non commutative.

Dans le cas envisag\'{e} au \S \ref{par-syst-classiques}, on a $\mathbf{D}=%
\mathbf{K}\left[ \partial \right] $ o\`{u} $\mathbf{K}$ est l'anneau $%
\mathcal{O}\left( 
%TCIMACRO{\U{211d} }%
%BeginExpansion
\mathbb{R}
%EndExpansion
\right) $\ des fonctions analytiques complexes sur $%
%TCIMACRO{\U{211d} }%
%BeginExpansion
\mathbb{R}
%EndExpansion
$, ou un sous-anneau de celui-ci. On peut aussi consid\'{e}rer des "syst\`{e}%
mes multidimensionnels" en prenant $\mathbf{D}=\mathbf{K}\left[ \partial
_{1},...,\partial _{n}\right] ,$ $\partial _{i}\triangleq \partial /\partial
x_{i},$ o\`{u} $\mathbf{K}$ est un sous-anneau de l'anneau $\mathcal{O}%
\left( \Omega \right) $\ des fonctions analytiques sur un ouvert $\Omega $ de%
$\ 
%TCIMACRO{\U{211d} }%
%BeginExpansion
\mathbb{R}
%EndExpansion
^{n}.$ On peut enfin supposer que $\mathbf{D}$ est une alg\`{e}bre de
convolution, etc.

Supposons pour fixer les id\'{e}es que $\mathbf{D}=\mathbf{K}\left[ \partial %
\right] $ o\`{u} $\mathbf{K=P}$, l'anneau des polyn\^{o}mes par rapport \`{a}
la variable $t$ (qui d\'{e}signe le temps), \`{a} coefficients complexes.
Remarquons tout d'abord que l'anneau $\mathbf{D}$ est non commutatif. En
effet, d'apr\`{e}s le r\`{e}gle de Leibniz,%
\begin{equation*}
\partial \left( tw\right) =w+t\partial w
\end{equation*}%
et comme ceci est valable pour toute variable $w$ on obtient%
\begin{equation}
\partial t-t\partial =1.  \label{regle-commutation}
\end{equation}%
Cet anneau $\mathbf{D}$ est isomorphe \`{a} la premi\`{e}re alg\`{e}bre de
Weyl habituellement not\'{e}e $A_{1}\left( 
%TCIMACRO{\U{2102} }%
%BeginExpansion
\mathbb{C}
%EndExpansion
\right) .$

Soit $W$ un $\mathbf{D}$-module \`{a} gauche et le "behavior"%
\begin{equation}
\mathfrak{B}_{S}\left( W\right) =\left\{ \mathbf{w}\in W^{k}:R\mathbf{w}%
=0\right\} .  \label{eq-behavior}
\end{equation}%
Dans le m\^{e}me esprit qu'au \S \ref{par-revol-groth},%
\begin{equation*}
\mathfrak{B}_{S}:W\mapsto \mathfrak{B}_{S}\left( W\right)
\end{equation*}%
est un foncteur contravariant de la cat\'{e}gorie des $\mathbf{D}$-modules 
\`{a} gauche dans la cat\'{e}gorie $\mathbf{Vect}_{%
%TCIMACRO{\U{211d} }%
%BeginExpansion
\mathbb{R}
%EndExpansion
}$ des espaces vectoriels sur $%
%TCIMACRO{\U{211d} }%
%BeginExpansion
\mathbb{R}
%EndExpansion
.$ Cette observation conduit \`{a} envisager l'\'{e}tude de $\mathfrak{B}%
_{S}\left( W\right) $ en deux \'{e}tapes:

\begin{itemize}
\item[(i)] l'\'{e}tude du foncteur $\mathfrak{B}_{S};$

\item[(ii)] l'\'{e}tude de "l'immersion" particuli\`{e}re%
\begin{equation*}
\mathfrak{B}_{S}\left( W\right) \rightarrow W^{k}.
\end{equation*}
\end{itemize}

L'\'{e}tape (i) rel\`{e}ve exclusivement de l'Alg\`{e}bre, tandis que l'\'{e}%
tape (ii) fait appel aux outils de l'Analyse, ce qui explique que cette
approche porte le nom d'\emph{Analyse Alg\'{e}brique}.

\section{$\mathbf{D}$-modules}

Il suffit maintenant de reprendre la br\`{e}ve pr\'{e}sentation qui a \'{e}t%
\'{e} faite plus haut du point de d\'{e}part de la G\'{e}om\'{e}trie Alg\'{e}%
brique moderne, en changeant l\'{e}g\`{e}rement de langage. Les \'{e}l\'{e}%
ments de $\mathfrak{B}_{S}\left( W\right) $ sont les \'{e}l\'{e}ments de $%
W^{k}$ qui s'annulent sur 
\begin{equation*}
\limfunc{im}\nolimits_{\mathbf{D}}\left( \bullet R\right) \triangleq \mathbf{%
D}^{1\times q}R,
\end{equation*}%
($\bullet R$ d\'{e}signant la multiplication \`{a} droite par $R$),
autrement dit%
\begin{equation*}
\mathfrak{B}_{S}\left( W\right) =\left\{ \mathbf{w}\in W^{k}:r\,\mathbf{w}%
=0,\forall r\in \limfunc{im}\nolimits_{\mathbf{D}}\left( \bullet R\right)
.\right\}
\end{equation*}%
Il existe un isomorphisme fonctoriel%
\begin{equation*}
\mathfrak{B}_{S}\left( W\right) \cong \limfunc{Hom}\nolimits_{\mathbf{D}%
}\left( \mathbf{D}^{1\times k}/\limfunc{im}\nolimits_{\mathbf{D}}\left(
\bullet R\right) ,W\right) .
\end{equation*}%
Soit%
\begin{equation*}
M=\mathbf{D}^{1\times k}/\limfunc{im}\nolimits_{\mathbf{D}}\left( \bullet
R\right) =\limfunc{coker}\nolimits_{\mathbf{D}}\left( \bullet R\right) .
\end{equation*}%
On peut maintenant identifier canoniquement $\mathfrak{B}_{S}\left( W\right) 
$ et $\limfunc{Hom}\nolimits_{\mathbf{D}}\left( M,W\right) .$

Tout $\mathbf{D}$-module \`{a} gauche de pr\'{e}sention finie est de la
forme $M=\limfunc{coker}\nolimits_{\mathbf{D}}\left( \bullet R\right) $, o%
\`{u} $R$ est une matrice (finie) \`{a} coefficients dans $\mathbf{D}.$ Les $%
\mathfrak{B}_{S}$ s'identifient aux foncteurs repr\'{e}sentables%
\begin{equation*}
\mathfrak{B}_{M}:W\mapsto \limfunc{Hom}\nolimits_{\mathbf{D}}\left(
M,W\right) .
\end{equation*}

Le foncteur $M\mapsto \mathfrak{B}_{M}$ est contravariant et la cat\'{e}%
gorie des foncteurs%
\begin{equation*}
_{\mathbf{D}}\mathbf{Mod}^{fp}\rightarrow \mathbf{Vect}_{%
%TCIMACRO{\U{211d} }%
%BeginExpansion
\mathbb{R}
%EndExpansion
}
\end{equation*}%
associ\'{e}s aux "behaviors" $\left( \ref{eq-behavior}\right) $ est \'{e}%
quivalente \`{a} la cat\'{e}gorie oppos\'{e}e de la cat\'{e}gorie des $%
\mathbf{D}$-modules \`{a} gauche de pr\'{e}sentation finie $_{\mathbf{D}}%
\mathbf{Mod}^{fp}$ en associant \`{a} tout $\mathbf{D}$-module $M$ de pr\'{e}%
sentation finie le foncteur $\mathfrak{B}_{M}.$

Nous sommes maintenant amen\'{e}s \`{a} identifier un syst\`{e}me lin\'{e}%
aire (sur l'anneau $\mathbf{D}$) et le $\mathbf{D}$-module (de pr\'{e}%
sentation finie) associ\'{e} \cite{Kashiwara}, \cite{Fliess-90} et \`{a}
adopter le point de vue suivant:

\emph{Le but initial de la Th\'{e}orie des Syst\`{e}mes moderne \'{e}quivaut 
\`{a} l'\'{e}tude des }$D$\emph{-modules de repr\'{e}sentation finie }$M$.

Par exemple, en supposant que $\mathbf{D}$ est un anneau d'Ore, on est
conduit \`{a} la d\'{e}finition suivante: le syst\`{e}me $S$ est commandable
si le $\mathbf{D}$-module associ\'{e} $M$ est sans torsion \cite{Fliess-90}, 
\cite{Pillai-Shankar}.

\subsection{Dualit\'{e}}

Le $\mathbf{D}$-module de pr\'{e}sentation finie $M$ constitue une repr\'{e}%
sentation intrins\`{e}que des \'{e}qutions du syst\`{e}me, tandis que $%
\mathfrak{B}_{M}\left( W\right) $ est l'ensemble de toutes les solutions
possibles, dans $W^{k},$ \`{a} ces \'{e}quations. Le risque est que la
structure alg\'{e}brique contenue dans $M$ soit beaucoup plus riche que la
structure des solutions, auquel cas l'\'{e}tude alg\'{e}brique de $M$ ne
fournirait que des r\'{e}sultats n'ayant aucune signification concr\`{e}te.

Une question essentielle est donc de savoir si la connaissance de $\mathfrak{%
B}_{M}\left( W\right) $ entra\^{\i}ne celle de $M$, autrement dit si le
foncteur $M\mapsto \mathfrak{B}_{M}\left( W\right) $ est injectif. Il suffit
pour cela que ce foncteur soit \emph{fid\`{e}le}, autrement dit que $W$ soit
cog\'{e}n\'{e}rateur dans la cat\'{e}gorie $_{\mathbf{D}}\mathbf{Mod}$ des $%
\mathbf{D}$-modules \`{a} gauche \cite{Oberst-multidimensional}. Une
condition suffisante moins restrictive est que $W\in {}_{\mathbf{D}}\mathbf{%
Mod}$ soit cog\'{e}n\'{e}rateur pour la sous-cat\'{e}gorie pleine $_{\mathbf{%
D}}\mathbf{Mod}^{fp}$ de $_{\mathbf{D}}\mathbf{Mod}$, ce qui conduit \`{a} d%
\'{e}finir dans une cat\'{e}gorie semi-ab\'{e}lienne $\mathcal{C}$ la notion
de cog\'{e}n\'{e}rateur pour une sous-cat\'{e}gorie $\mathcal{D}$ de $%
\mathcal{C}$ \cite{Bourles-Oberst}.

Il est \'{e}vident, par exemple, qu'en choisissant $W=0,$ on obtient $%
\mathfrak{B}_{M}\left( W\right) =0$ quel que soit le $\mathbf{D}$-module $M.$
Des cas moins triviaux conduisent \`{a} des questions d'Analyse fine. Par
exemple, consid\'{e}rons un anneau semblable \`{a} la premi\`{e}re alg\`{e}%
bre de Weyl $A_{1}\left( 
%TCIMACRO{\U{2102} }%
%BeginExpansion
\mathbb{C}
%EndExpansion
\right) $ mais o\`{u} l'anneau des coefficients $\mathbf{P}=%
%TCIMACRO{\U{2102} }%
%BeginExpansion
\mathbb{C}
%EndExpansion
\left[ t\right] $ est remplac\'{e} par $\mathcal{R}\left( \Omega \right) =%
%TCIMACRO{\U{2102} }%
%BeginExpansion
\mathbb{C}
%EndExpansion
\left( t\right) \cap \mathcal{O}\left( \Omega \right) ,$ o\`{u} $\Omega $
est un intervalle ouvert non vide de la droite r\'{e}elle. Notons $%
A_{0}\left( \Omega \right) $ l'anneau d'op\'{e}rateurs diff\'{e}rentiels $%
\mathcal{R}\left( \Omega \right) \left[ \partial \right] ,$ muni de la r\`{e}%
gle de commutation $\left( \ref{regle-commutation}\right) .$ On v\'{e}rifie
que $A_{0}\left( \Omega \right) $ est, comme $A_{1}\left( 
%TCIMACRO{\U{2102} }%
%BeginExpansion
\mathbb{C}
%EndExpansion
\right) ,$ un anneau de Dedekind simple non commutatif. Consid\'{e}rons le
syst\`{e}me%
\begin{equation}
\left( t^{3}\partial +2\right) w=0.  \label{eq-Schwartz}
\end{equation}%
o\`{u} $t^{3}\partial +2\in \mathbf{D}\triangleq A_{0}\left( 
%TCIMACRO{\U{2102} }%
%BeginExpansion
\mathbb{C}
%EndExpansion
\right) .$ \ Le $\mathbf{D}$-module%
\begin{equation*}
M=\mathbf{D}/\mathbf{D}\left( t^{3}\partial +2\right)
\end{equation*}%
est un module de torsion non r\'{e}duit \`{a} $0$, il s'agit donc d'un syst%
\`{e}me non commandable. Assez naturellement, on est tent\'{e} de chercher
les solutions dans $W=\mathcal{D}^{\prime }\left( 
%TCIMACRO{\U{211d} }%
%BeginExpansion
\mathbb{R}
%EndExpansion
\right) ,$ l'espace des distributions sur la droite r\'{e}elle. Mais l'\'{e}%
quation $\left( \ref{eq-Schwartz}\right) $ admet $0$ comme unique solution
dans $\mathcal{D}^{\prime }\left( 
%TCIMACRO{\U{211d} }%
%BeginExpansion
\mathbb{R}
%EndExpansion
\right) $ (\cite{Schwartz-distrib}, ((V,6;15)). Par cons\'{e}quent l'espace
des distributions ne fournit pas une dualit\'{e} syst\`{e}mes $%
\longleftrightarrow $ "behaviors" lorsque $\mathbf{D}=A_{0}\left( 
%TCIMACRO{\U{2102} }%
%BeginExpansion
\mathbb{C}
%EndExpansion
\right) $.

En revanche, l'espace $B\left( 
%TCIMACRO{\U{211d} }%
%BeginExpansion
\mathbb{R}
%EndExpansion
\right) $ des hyperfonctions de Sato \cite{Sato} est un $A_{0}\left( 
%TCIMACRO{\U{2102} }%
%BeginExpansion
\mathbb{C}
%EndExpansion
\right) $-module cog\'{e}n\'{e}rateur \cite{Frohler-Oberst} (et, plus pr\'{e}%
cis\'{e}ment, est un "large injective cogenerator") et permet donc d'obtenir
la dualit\'{e} recherch\'{e}e.

Voici un aper\c{c}u de l'approche moderne de la Th\'{e}orie des Syst\`{e}%
mes, qui fait l'objet de l'ouvrage \cite{Bourles-Marinescu}, et de mes
recherches actuelles.

\end{document}